\documentclass{amsart}%
\usepackage{amssymb}
\usepackage{amsfonts}
\usepackage{amsmath}
\usepackage{graphicx}%
\newtheorem{theorem}{Theorem}
\theoremstyle{plain}

\newtheorem{notation}{Notation}

\numberwithin{equation}{section}
\begin{document}
\title{Non-Archimedean Models of Morphogenesis}
\author{W. A. Z\'{u}\~{n}iga-Galindo}
\address{University of Texas Rio Grande Valley\\
School of Mathematical \& Statistical Sciences\\
One West University Blvd\\
Brownsville, TX 78520, United States and Centro de Investigaci\'{o}n y de
Estudios Avanzados del Instituto Polit\'{e}cnico Nacional\\
Departamento de Matem\'{a}ticas, Unidad Quer\'{e}taro\\
Libramiento Norponiente \#2000, Fracc. Real de Juriquilla. Santiago de
Quer\'{e}taro, Qro. 76230\\
M\'{e}xico.}
\email{wilson.zunigagalindo@utrgv.edu, wazuniga@math.cinvestav.edu.mx}
\thanks{The author was partially supported by Conacyt Grant No. 217367 (Mexico), and
by the Debnath Endowed Professorship (UTRGV, USA)}

\begin{abstract}
We study a $p$-adic reaction-diffusion system and the associated Turing
patterns. We establish an instability criteria and show that the Turing
patterns are not classical patterns consisting of alternating domains. Instead
of this, a Turing pattern consists of several domains (clusters), each of them
supporting a different pattern but with the same parameter values. This type
of patterns are typically produced by reaction-diffusion equations on large networks.

\end{abstract}
\maketitle

\section{Introduction}

In 1952, A. Turing proposed that under certain conditions chemicals can react
and diffuse in such a way as to produce steady state heterogeneous spatial
patterns of chemical (or morphogen) concentration. In the case of two
morphogens interacting, the model proposed by Turing has the form:%
\begin{equation}
\left\{
\begin{array}
[c]{cc}%
\frac{\partial u}{\partial t}(x,t)= & \gamma f\left(  u,v\right)
+\frac{\partial^{2}u}{\partial x^{2}}(x,t)\\
& \\
\frac{\partial v}{\partial t}(x,t)= & \gamma g\left(  u,v\right)
+d\frac{\partial^{2}v}{\partial x^{2}}(x,t),
\end{array}
\right.  \label{EQ_0}%
\end{equation}
with suitable boundary conditions. In this model \ $u$, $v$ represent the
concentrations\ of the two morphogens, $f$ and $g$ represent the reaction
kinetics, and $d$ is \ the ratio of diffusion, and $\gamma$ represents the
relative strength of the reaction terms.

Typically $u\left(  x,t\right)  $ and $v(x,t)$ are interpreted to be the local
densities of the activator and inhibitor species. Functions $f\left(
u,v\right)  $ and $g\left(  u,v\right)  $ specify the local dynamics of the
activator, which autocatalytically enhances its own production, and of the
inhibitor, which suppresses the activator growth. The Turing instability
occurs when the parameter $d$ exceeds a threshold $d_{c}$, \cite{Turing},
\cite{Murray-II}. This event drives to a spontaneous development of a spatial
pattern formed by alternating activator-rich and activator-poor patches.
Turing instability in activator-inhibitor systems establishes a paradigm of
non-equilibrium self-organization, which has been extensively studied for
biological and chemical processes.

In the 70s, Othmer and Scriven started the study of the Turing instability in
network-organized systems \cite{Othmer et al.}-\cite{Othmer et al 2}. Since
then, reaction-diffusion models on networks have been studied intensively, see
e.g. \cite{Ambrosio et al}, \cite{Boccaletti et al}, \cite{Chung}, \cite{Ide},
\cite{Mocarlo}, \cite{Mugnolo}, \cite{Nakao-Mikhalov}, \cite{Othmer et al.},
\cite{Othmer et al 2}, \cite{Slavova et al}, \cite{Van Mieghem}, \cite{von
Below}, \cite{Zhao}, \cite{Zuniga-JMAA} and the references therein.

In the discrete case, the continuous media is replaced by a network (an
unoriented graph $\mathcal{G}$, which plays the role of discrete media). The
analog of operator $\Delta$ is the Laplacian of the graph $\mathcal{G}$, which
is defined as%
\begin{equation}
\left[  L_{JI}\right]  _{J,I\in V(\mathcal{G})}=\left[  A_{JI}-\gamma
_{I}\delta_{JI}\right]  _{J,I\in V(\mathcal{G})}, \label{EQ_Matrix_L}%
\end{equation}
where $\left[  A_{JI}\right]  _{J,I\in V(\mathcal{G})}$ is the adjacency
matrix of $\mathcal{G}$ and $\gamma_{I}$ is the degree of $I$. The network
analogue of (\ref{EQ_0}) is
\begin{equation}
\left\{
\begin{array}
[c]{l}%
\frac{\partial u_{J}}{\partial t}=f(u_{J},v_{J})+\varepsilon%
{\textstyle\sum\limits_{I}}
L_{JI}u_{I}\\
\\
\frac{\partial v_{J}}{\partial t}=g(u_{J},v_{J})+\varepsilon d%
{\textstyle\sum\limits_{I}}
L_{JI}v_{I}.
\end{array}
\right.  \label{EQ_1}%
\end{equation}

In the last fifty years, Turing patterns produced by reaction-diffusion
systems on networks have been studied intensively, see e.g. \cite{Boccaletti
et al}, \cite{Ide}, \cite{Mocarlo}, \cite{Nakao-Mikhalov}, \cite{Othmer et
al.}, \cite{Othmer et al 2}, \cite{Zhao}, \cite{Zuniga-JMAA} \ and the
references therein. Nowadays, there is a large amount of experimental results
about the behavior of these systems, obtained mainly via computer simulations
using large random networks. The investigations of the Turing patterns for
large random networks have revealed that, whereas the Turing criteria remain
essentially the same, as in the classical case, the properties of the emergent
patterns are very different. In \cite{Nakao-Mikhalov}, by using a physical
argument, Nakao and Mikhailov establish\ that Turing patterns with alternating
domains cannot exist in the network case, and only several domains (clusters)
occur. Multistability, that is, coexistence of a number of different patterns
with the same parameter values, is typically found and hysteresis phenomena
are observed. They used mean-field approximation to understand the Turing
patterns when $d>d_{c}$, and proposed that the mean-field approximation is the
natural framework to understand the peculiar behavior of the Turing patterns
on networks. This program was carried out in \cite{Zuniga-JMAA} using $p$-adic analysis.

In \cite{Zuniga-JMAA} we established that the $p$-adic reaction-diffusion
system%
\begin{equation}
\left\{
\begin{array}
[c]{l}%
\frac{\partial u\left(  x,t\right)  }{\partial t}=f(u\left(  x,t\right)
,v\left(  x,t\right)  )+\varepsilon%
{\textstyle\int\limits_{\mathcal{K}_{N}}}
\left(  u\left(  y,t\right)  -u\left(  x,t\right)  \right)  J_{N}(x,y)dy\\
\\
\frac{\partial v\left(  x,t\right)  }{\partial t}=g(u\left(  x,t\right)
,v\left(  x,t\right)  )+\varepsilon d%
{\textstyle\int\limits_{\mathcal{K}_{N}}}
\left(  u\left(  y,t\right)  -u\left(  x,t\right)  \right)  J_{N}(x,y)dy,
\end{array}
\right.  \label{EQ_1A}%
\end{equation}
where $\mathcal{K}_{N}$ (an open compact subset) and $J_{N}(x,y)$ depends on
$\left[  A_{JI}\right]  _{J,I\in V(\mathcal{G})}$, is a good $p$-adic
continuous approximation of (\ref{EQ_1}). The Turing instability criteria
\ for (\ref{EQ_1A}) is essentially the same as in the classical case, but the
qualitative description of the Turing patterns is the one given by Nakao and
Mikhailov \ in terms of clustering and multistability.

In this article we study the following $p$-adic reaction-diffusion system:%
\begin{equation}
\left\{
\begin{array}
[c]{l}%
\frac{\partial u}{\partial t}(x,t)=\gamma f\left(  u,v\right)  -\boldsymbol{D}%
_{x}^{\alpha}u(x,t)\text{;}\\
\\
\frac{\partial v}{\partial t}(x,t)=\gamma g\left(  u,v\right)
-d\boldsymbol{D}_{x}^{\alpha}v(x,t)\text{,}%
\end{array}
\right.  \label{EQ_1B}%
\end{equation}
where $x\in\mathbb{Q}_{p}$, $t\geq0$, and $\boldsymbol{D}_{x}^{\alpha}$ is the
Vladimirov operator. Since want to study self-organization patterns we need a
`zero flux boundary condition', which is
\[
u(x,t),v(x,t)\equiv0\text{ for }x\in\mathbb{Q}_{p}\smallsetminus B_{M}\text{,
for any }t\geq0,
\]
where $B_{M}$ is the ball of radius $p^{M}$ around the origin. We establish a
Turing instability criteria \ for (\ref{EQ_1B}) and show that the Turing
pattern can be described as a couple of convergent series involving the
functions%
\begin{align*}
&  e^{\lambda t}p^{\frac{-r}{2}}\cos\left(  \left\{  p^{-1}j\left(
p^{r}x-n\right)  \right\}  _{p}\right)  \Omega\left(  \left\vert
p^{r}x-n\right\vert _{p}\right)  ,\\
&  \text{ }p^{\frac{-r}{2}}\sin\left(  \left\{  p^{-1}j\left(  p^{r}%
x-n\right)  \right\}  _{p}\right)  \Omega\left(  \left\vert p^{r}%
x-n\right\vert _{p}\right)  ,\text{ }%
\end{align*}
see Theorem \ref{Theorem1}. The Turing patterns attached to (\ref{EQ_1B}) are
not classical patterns consisting of alternating domains. Instead of this, a
Turing pattern consists of several domains (clusters), each of them supporting
a different pattern but with the same parameter values (multistability). It is
important to mention that the spectra of $\boldsymbol{D}_{x}^{\alpha}$
consists of a sequence of positive eigenvalues with infinite multiplicity,
while the spectra of the operators considered in \cite[Theorem 10.1]%
{Zuniga-JMAA} consists of a finite number of non-negative eigenvalues with
finite multiplicities. In Section \ref{Section_Discrete_models} we \ construct
discretization of system (\ref{EQ_1B}) of type (\ref{EQ_1}). But, In this case
the matrix of the \ discrete is not related to adjacency matrix of a graph,
see (\ref{EQ_A_matrix}).

\section{\label{Fourier Analysis}$p$-Adic Analysis: Essential Ideas}

In this section we collect some basic results about $p$-adic analysis that
will be used in the article. For an in-depth review of the $p$-adic analysis
the reader may consult \cite{Alberio et al}, \cite{Taibleson}, \cite{V-V-Z}.

\subsection{The field of $p$-adic numbers}

Along this article $p$ will denote a prime number. The field of $p-$adic
numbers $%
\mathbb{Q}
_{p}$ is defined as the completion of the field of rational numbers
$\mathbb{Q}$ with respect to the $p-$adic norm $|\cdot|_{p}$, which is defined
as
\[
\left\vert x\right\vert _{p}=\left\{
\begin{array}
[c]{lll}%
0 & \text{if} & x=0\\
&  & \\
p^{-\gamma} & \text{if} & x=p^{\gamma}\frac{a}{b}\text{,}%
\end{array}
\right.
\]
where $a$ and $b$ are integers coprime with $p$. The integer $\gamma:=ord(x)
$, with $ord(0):=+\infty$, is called the\textit{\ }$p-$\textit{adic order of}
$x$.

Any $p-$adic number $x\neq0$ has a unique expansion of the form
\[
x=p^{ord(x)}\sum_{j=0}^{\infty}x_{j}p^{j},
\]
where $x_{j}\in\{0,\dots,p-1\}$ and $x_{0}\neq0$. By using this expansion, we
define \textit{the fractional part of }$x\in\mathbb{Q}_{p}$, denoted
$\{x\}_{p}$, as the rational number
\[
\left\{  x\right\}  _{p}=\left\{
\begin{array}
[c]{lll}%
0 & \text{if} & x=0\text{ or }ord(x)\geq0\\
&  & \\
p^{ord(x)}\sum_{j=0}^{-ord_{p}(x)-1}x_{j}p^{j} & \text{if} & ord(x)<0.
\end{array}
\right.
\]
In addition, any non-zero $p-$adic number can be represented uniquely as
$x=p^{ord(x)}ac\left(  x\right)  $ where $ac\left(  x\right)  =\sum
_{j=0}^{\infty}x_{j}p^{j}$, $x_{0}\neq0$, is called the \textit{angular
component} of $x$. Notice that $\left\vert ac\left(  x\right)  \right\vert
_{p}=1$.

For $r\in\mathbb{Z}$, denote by $B_{r}(a)=\{x\in%
\mathbb{Q}
_{p};\left\vert x-a\right\vert _{p}\leq p^{r}\}$ \textit{the ball of radius
}$p^{r}$ \textit{with center at} $a\in%
\mathbb{Q}
_{p}$, and take $B_{r}(0):=B_{r}$. The ball $B_{0}$ equals $\mathbb{Z}_{p}$,
\textit{the ring of }$p-$\textit{adic integers of }$%
\mathbb{Q}
_{p}$. We also denote by $S_{r}(a)=\{x\in\mathbb{Q}_{p};|x-a|_{p}=p^{r}\}$
\textit{the sphere of radius }$p^{r}$ \textit{with center at} $a\in%
\mathbb{Q}
_{p}$, and take $S_{r}(0):=S_{r}$. We notice that $S_{0}^{1}=\mathbb{Z}%
_{p}^{\times}$ (the group of units of $\mathbb{Z}_{p}$). The balls and spheres
are both open and closed subsets in $%
\mathbb{Q}
_{p}$. In addition, two balls in $%
\mathbb{Q}
_{p}$ are either disjoint or one is contained in the other.

The metric space $\left(
\mathbb{Q}
_{p},\left\vert \cdot\right\vert _{p}\right)  $ is a complete ultrametric
space. As a topological space $\left(
\mathbb{Q}
_{p},|\cdot|_{p}\right)  $ is totally disconnected, i.e. the only connected
\ subsets of $%
\mathbb{Q}
_{p}$ are the empty set and the points. In addition, $\mathbb{Q}_{p}$\ is
homeomorphic to a Cantor-like subset of the real line, see e.g. \cite{Alberio
et al}, \cite{V-V-Z}. A subset of $\mathbb{Q}_{p}$ is compact if and only if
it is closed and bounded in $\mathbb{Q}_{p}$, see e.g. \cite[Section
1.3]{V-V-Z}, or \cite[Section 1.8]{Alberio et al}. The balls and spheres are
compact subsets. Thus $\left(
\mathbb{Q}
_{p},|\cdot|_{p}\right)  $ is a locally compact topological space.

\begin{notation}
We use $\Omega\left(  p^{-r}|x-a|_{p}\right)  $ to denote the characteristic
function of the ball $B_{r}(a)=a+p^{-r}\mathbb{Z}_{p}$. For more general sets,
we use the notation $1_{A}$ for the characteristic function of a $A$.
\end{notation}

\subsection{Some function spaces}

A complex-valued function $\varphi$ defined on $%
\mathbb{Q}
_{p}$ is \textit{called locally constant} if for any $x\in%
\mathbb{Q}
_{p}$ there exist an integer $l(x)\in\mathbb{Z}$ such that
\begin{equation}
\varphi(x+x^{\prime})=\varphi(x)\text{ for }x^{\prime}\in B_{l(x)}.
\label{local_constancy_parameter}%
\end{equation}
\ A function $\varphi:%
\mathbb{Q}
_{p}\rightarrow\mathbb{C}$ is called a \textit{Bruhat-Schwartz function (or a
test function)} if it is locally constant with compact support. In this case,
we can take $l=l(\varphi)$ in (\ref{local_constancy_parameter}) independent of
$x$, the largest of such integers is called \textit{the parameter of local
constancy} of $\varphi$. The $\mathbb{C}$-vector space of Bruhat-Schwartz
functions is denoted by $\mathcal{D}:=\mathcal{D}(%
\mathbb{Q}
_{p})$. We denote by $\mathcal{D}_{\mathbb{R}}:=\mathcal{D}_{\mathbb{R}}(%
\mathbb{Q}
_{p}^{n})$, the $\mathbb{R}$-vector space of test functions.

Since $(\mathbb{Q}_{p},+)$ is a locally compact topological group, there
exists a Borel measure $dx$, called the Haar measure of $(\mathbb{Q}_{p},+)$,
unique up to multiplication by a positive constant. Furthermore, $\int
_{U}dx>0$ for every non-empty open set $U\subset\mathbb{Q}_{p}$, and
$\int_{E+z}dx=\int_{E}dx$ for every Borel set $E\subset\mathbb{Q}_{p}$, see
e.g. \cite[Chapter XI]{Halmos}. If we normalize this measure by the condition
$\int_{\mathbb{Z}_{p}}dx=1$, then $dx$ is unique. From now on we denote by
$dx$ the normalized Haar measure of $(\mathbb{Q}_{p},+)$.

Given $\rho\in\lbrack0,\infty)$, we denote by $L^{\rho}:=L^{\rho}\left(
\mathbb{Q}_{p}\right)  :=L^{\rho}\left(
\mathbb{Q}
_{p},dx\right)  ,$ the $%
\mathbb{C}
-$vector space of all the complex valued functions $g$ satisfying $\int_{%
\mathbb{Q}
_{p}}\left\vert g\left(  x\right)  \right\vert ^{\rho}dx<\infty$, and
$L^{\infty}\allowbreak:=L^{\infty}\left(
\mathbb{Q}
_{p}\right)  =L^{\infty}\left(
\mathbb{Q}
_{p},dx\right)  $ denotes the $%
\mathbb{C}
-$vector space of all the complex valued functions $g$ such that the essential
supremum of $|g|$ is bounded. The corresponding $\mathbb{R}$-vector spaces are
denoted as $L_{\mathbb{R}}^{\rho}\allowbreak:=L_{\mathbb{R}}^{\rho}\left(
\mathbb{Q}_{p}\right)  =L_{\mathbb{R}}^{\rho}\left(  \mathbb{Q}_{p},dx\right)
$, $1\leq\rho\leq\infty$.

\subsection{Fourier transform}

Set $\chi_{p}(y)=\exp(2\pi i\{y\}_{p})$ for $y\in%
\mathbb{Q}
_{p}$. The map $\chi_{p}(\cdot)$ is an additive character on $%
\mathbb{Q}
_{p}$, i.e. a continuous map from $\left(
\mathbb{Q}
_{p},+\right)  $ into $S$ (the unit circle considered as multiplicative group)
satisfying $\chi_{p}(x_{0}+x_{1})=\chi_{p}(x_{0})\chi_{p}(x_{1})$,
$x_{0},x_{1}\in%
\mathbb{Q}
_{p}$. The additive characters of $%
\mathbb{Q}
_{p}$ form an Abelian group which is isomorphic to $\left(
\mathbb{Q}
_{p},+\right)  $. The isomorphism is given by $\xi\rightarrow\chi_{p}(\xi x)$,
see e.g. \cite[Section 2.3]{Alberio et al}.

If $f\in L^{1}$ its Fourier transform is defined by
\[
(\mathcal{F}f)(\xi)=\int_{%
\mathbb{Q}
_{p}}\chi_{p}(\xi x)f(x)dx,\quad\text{for }\xi\in%
\mathbb{Q}
_{p}.
\]
We will also use the notation $\mathcal{F}_{x\rightarrow\xi}f$ and
$\widehat{f}$\ for the Fourier transform of $f$. The Fourier transform is a
linear isomorphism from $\mathcal{D}$ onto itself satisfying
\begin{equation}
(\mathcal{F}(\mathcal{F}f))(\xi)=f(-\xi), \label{FF(f)}%
\end{equation}
for every $f\in\mathcal{D},$ see e.g. \cite[Section 4.8]{Alberio et al}. If
$f\in L^{2}$, its Fourier transform is defined as
\[
(\mathcal{F}f)(\xi)=\lim_{k\rightarrow\infty}\int_{|x|_{p}\leq p^{k}}\chi
_{p}(\xi\cdot x)f(x)d^{n}x,\quad\text{for }\xi\in%
\mathbb{Q}
_{p}\text{,}%
\]
where the limit is taken in $L^{2}$. We recall that the Fourier transform is
unitary on $L^{2},$ i.e. $||f||_{L^{2}}=||\mathcal{F}f||_{L^{2}}$ for $f\in
L^{2}$ and that (\ref{FF(f)}) is also valid in $L^{2}$, see e.g. \cite[Chapter
$III$, Section 2]{Taibleson}.

\subsection{The Vladimirov Operator}

The Vladimirov pseudodifferential operator $\boldsymbol{D}^{\alpha}$,
$\alpha>0$, is defined as%
\begin{equation}
\boldsymbol{D}^{\alpha}\varphi\left(  x\right)  =\frac{1-p^{\alpha}%
}{1-p^{-\alpha-1}}\int_{\mathbb{Q}_{p}}|y|_{p}^{-\alpha-1}(\varphi
(x-y)-\varphi(x))\,dy\text{, for }\varphi\in\mathcal{D}. \label{8}%
\end{equation}
The right-hand side of (\ref{8}) makes sense for a wider class of functions,
for example, for locally constant functions $\varphi$ satisfying
\[
\int_{|x|_{p}\geq1}|x|_{p}^{-\alpha-d}|\varphi(x)|\,dx<\infty.
\]
Consequently, the constant functions are contained in the domain of
$\boldsymbol{D}^{\alpha}$, and that $\boldsymbol{D}^{\alpha}\varphi=0$, for
any constant function $\varphi$. On other hand,
\begin{equation}
\boldsymbol{D}^{\alpha}\varphi(x)=\mathcal{F}_{\xi\rightarrow x}^{-1}\left(
\left\vert \xi\right\vert _{p}^{\alpha}\mathcal{F}_{x\rightarrow\xi}%
\varphi\right)  \text{, for }\varphi\in\mathcal{D}\text{.}
\label{Taibleson_operator}%
\end{equation}
Finally in case in which the Vladimirov acts on functions depending on two
variables, $(x,t)\in\mathbb{Q}_{p}\times\mathbb{R}_{+}$, we will sue the
notation $\boldsymbol{D}_{x}^{\alpha}u(x,t)$ instead of $\boldsymbol{D}%
^{\alpha}u(x,t)$.

\subsubsection{\label{Section_spectrum}The spectrum of the operator
$\boldsymbol{D}^{\alpha}$}

The set of functions $\left\{  \Psi_{rnj}\right\}  $ defined as%
\begin{equation}
\Psi_{rnj}\left(  x\right)  =p^{\frac{-r}{2}}\chi_{p}\left(  p^{-1}j\left(
p^{r}x-n\right)  \right)  \Omega\left(  \left\vert p^{r}x-n\right\vert
_{p}\right)  ,\label{eq4}%
\end{equation}
where $r\in\mathbb{Z}$, $j\in\left\{  1,\cdots,p-1\right\}  $, and $n$ runs
through a fixed set of representatives of $\mathbb{Q}_{p}/\mathbb{Z}_{p}$, is
an orthonormal basis of $L^{2}(\mathbb{Q}_{p})$ consisting of eigenvectors of
operator $\boldsymbol{D}^{\alpha}$:%
\begin{equation}
\boldsymbol{D}^{\alpha}\Psi_{rnj}=p^{\left(  1-r\right)  \alpha}\Psi
_{rnj}\text{ for any }r\text{, }n\text{, }j\text{,}\label{eq5}%
\end{equation}

We set $\mathcal{L}\left(  B_{M}\right)  $ to be the $\mathbb{C}$-vector space
generated by the functions $\Psi_{rnj}\left(  x\right)  $ with support in
$B_{M}$, which are \ exactly those satisfying
\begin{equation}
r\leq M,\text{ }n\in p^{r-M}\mathbb{Z}_{p}\cap\mathbb{Q}_{p}/\mathbb{Z}%
_{p},\text{ \ }j\in\left\{  1,\cdots,p-1\right\}  . \label{conditions}%
\end{equation}
Notice that $\mathcal{L}\left(  B_{M}\right)  $ is a closed subspace of
$L^{2}\left(  B_{M}\right)  $. Futhermore, all the functions $\Psi
_{rnj}\left(  x\right)  $ satisfying (\ref{conditions}) are orthogonal to the
characteristic function of $B_{M}$, i.e. $\int\nolimits_{B_{M}}\Psi
_{rnj}\left(  x\right)  dx=0$. Notice that $L^{2}\left(  B_{M}\right)
=\mathbb{C}\Omega\left(  p^{-M}\left\vert x\right\vert _{p}\right)
{\textstyle\bigoplus}
L_{0}^{2}\left(  B_{M}\right)  $, where
\[
L_{0}^{2}\left(  B_{M}\right)  =\left\{  f\in L^{2}\left(  B_{M}\right)
;\int_{B_{M}}fdx=0\right\}  .
\]

\subsection{\label{Section_Eigenvalue_problems}Two spectral problems}

Consider the spectral problem:%
\begin{equation}
\left\{
\begin{array}
[c]{ll}%
\boldsymbol{D}^{\alpha}\theta\left(  x\right)  =\kappa\theta\left(  x\right)
\text{,} & \kappa\in\mathbb{R}\\
& \\
\theta\in L_{\mathbb{R}}^{2}\left(  \mathbb{Q}_{p}\right)  . &
\end{array}
\right.  \label{E_Value_1}%
\end{equation}
The functions $\Psi_{rnj}\left(  x\right)  $ are complex-valued eigenfunctions
of (\ref{E_Value_1}) with eigenvalues $\kappa\in\left\{  p^{\left(
1-r\right)  \alpha};r\in\mathbb{Z}\right\}  $. Notice that each eigenvalues
has infinite multiplicity. Therefore
\begin{align}
&  p^{\frac{-r}{2}}\cos\left(  \left\{  p^{-1}j\left(  p^{r}x-n\right)
\right\}  _{p}\right)  \Omega\left(  \left\vert p^{r}x-n\right\vert
_{p}\right)  \text{, }\nonumber\\
&  p^{\frac{-r}{2}}\sin\left(  \left\{  p^{-1}j\left(  p^{r}x-n\right)
\right\}  _{p}\right)  \Omega\left(  \left\vert p^{r}x-n\right\vert
_{p}\right)  \text{,} \label{Eq_real_basis}%
\end{align}
with $r$, $j$, $n$ as before, are real-valued eigenfunctions of
(\ref{E_Value_1}) with $\kappa=p^{\left(  1-r\right)  \alpha}$. Notice that
the eigenfunctions are not completely determined by the `wavenumber' $\kappa$.
The functions of the type (\ref{Eq_real_basis}) form a basis of $L_{\mathbb{R}%
}^{2}\left(  \mathbb{Q}_{p}\right)  $ (which is not necessarily orthonormal).
More precisely, $f(x)=\sum_{rnj}A_{rnj}\Psi_{rnj}\left(  x\right)  \in
L_{\mathbb{R}}^{2}\left(  \mathbb{Q}_{p}\right)  $ admits an expansion of the
form%
\begin{align}
&  \sum\limits_{rnj}p^{\frac{-r}{2}}\operatorname{Re}(A_{rnj})\cos\left(
\left\{  p^{-1}j\left(  p^{r}x-n\right)  \right\}  _{p}\right)  \Omega\left(
\left\vert p^{r}x-n\right\vert _{p}\right)  -\label{E_expansio}\\
&  \sum\limits_{rnj}p^{\frac{-r}{2}}\operatorname{Im}(A_{rnj})\sin\left(
\left\{  p^{-1}j\left(  p^{r}x-n\right)  \right\}  _{p}\right)  \Omega\left(
\left\vert p^{r}x-n\right\vert _{p}\right)  ,\nonumber
\end{align}
where%
\begin{align*}
\operatorname{Re}(A_{rnj})  &  =p^{\frac{-r}{2}}\int\limits_{\mathbb{Q}_{p}%
}f\left(  x\right)  \cos\left(  \left\{  p^{-1}j\left(  p^{r}x-n\right)
\right\}  _{p}\right)  \Omega\left(  \left\vert p^{r}x-n\right\vert
_{p}\right)  dx\text{, }\\
\operatorname{Im}(A_{rnj})  &  =p^{\frac{-r}{2}}\int\limits_{\mathbb{Q}_{p}%
}f\left(  x\right)  \sin\left(  \left\{  p^{-1}j\left(  p^{r}x-n\right)
\right\}  _{p}\right)  \Omega\left(  \left\vert p^{r}x-n\right\vert
_{p}\right)  dx\text{. }%
\end{align*}

Now we consider the eigenvalue problem:%
\begin{equation}
\left\{
\begin{array}
[c]{ll}%
\boldsymbol{D}_{x}^{\alpha}\theta\left(  x\right)  =\kappa\theta\left(
x\right)  \text{,} & \kappa\in\mathbb{R}\\
& \\
\theta\in L_{\mathbb{R}}^{2}\left(  B_{M}\right)  \cap\mathcal{L}\left(
B_{M}\right)  , & M\in\mathbb{Z}\text{.}%
\end{array}
\right.  \label{E_Value_2}%
\end{equation}
All the functions $\operatorname{Re}\left(  \Psi_{rnj}\left(  x\right)
\right)  $ satisfying (\ref{conditions}) are solutions of (\ref{E_Value_2}),
with $\kappa=p^{\left(  1-r\right)  \alpha}$. Now, any function $f\in
L_{\mathbb{R}}^{2}\left(  B_{M}\right)  \cap\mathcal{L}\left(  B_{M}\right)  $
admits a Fourier expansion of the form
\begin{equation}
f(x)=\sum\nolimits_{rnj}C_{rnj}\Psi_{rnj}\left(  x\right)  =\sum
\nolimits_{rnj}\operatorname{Re}\left(  C_{rnj}\Psi_{rnj}\left(  x\right)
\right)  , \label{Eq_30}%
\end{equation}
where the $\Psi_{rnj}\left(  x\right)  $s run through all the wavelets
supported in the ball $B_{M}$. Therefore, (\ref{Eq_30}) is a solution of
(\ref{E_Value_2}).

Finally we notice that%
\begin{equation}
\boldsymbol{D}^{\alpha}\Omega\left(  p^{-M}\left\vert x\right\vert
_{p}\right)  =\left\{
\begin{array}
[c]{lll}%
\frac{\left(  1-p^{-1}\right)  p^{-\alpha M}}{1-p^{-\alpha-1}} & \text{if} &
\left\vert x\right\vert _{p}\leq p^{M}\\
&  & \\
\frac{\left(  1-p^{\alpha}\right)  p^{M}}{1-p^{-\alpha-1}}\frac{1}{\left\vert
x\right\vert _{p}^{\alpha+1}} &  & \left\vert x\right\vert _{p}>p^{M}%
\end{array}
\right.  \label{Eq_D_alpha}%
\end{equation}
implies that $\Omega\left(  p^{-M}\left\vert x\right\vert _{p}\right)  $ is
not a solution of (\ref{E_Value_2}).

\subsection{The $p$-adic heat equation}

The evolution equation
\begin{equation}
\frac{\partial u(x,t)}{\partial t}+(\boldsymbol{D}_{x}^{\alpha}u)(x,t)=0,\quad
x\in\mathbb{Q}_{p},\quad t\geq0, \label{Eq_1}%
\end{equation}
is the $p$-adic heat equation. The analogy with the classical heat equation
comes from the fact that the solution of the initial value problem attached to
(\ref{Eq_1}) with initial datum $u(x,0)=\varphi(x)\in\mathcal{D}_{\mathbb{R}}$
is given by%
\[
u(x,t)=\int_{\mathbb{Q}_{p}}Z(x-y,t)\varphi(x)\,dx,
\]
where%
\[
Z(x,t):=\int_{\mathbb{Q}_{p}}\chi_{p}(-x\xi)e^{-t|\xi|_{p}^{\alpha}}%
\,d\xi\text{ for }t>0\text{,}%
\]
is the $p$\textit{-adic heat kernel}. $Z\left(  x,t\right)  $ is a transition
density of a time and space homogeneous Markov process which is bounded, right
continuous and has no discontinuities other than jumps, see e.g. \cite[Theorem
16]{Zuniga-LNM-2016}.

We now review the classical case. For $\gamma>0$, consider the fractional
Laplacian:%
\[
\left(  -\Delta\right)  ^{\gamma}\varphi\left(  x\right)  =\mathcal{F}%
_{\xi\rightarrow x}^{-1}\left(  |\xi|_{\mathbb{R}}^{\gamma}\mathcal{F}%
_{x\rightarrow\xi}\varphi\right)  ,
\]
where $\varphi$ is a Schwartz function, $\mathcal{F}$ denotes the Fourier
transform in the group $\left(  \mathbb{R},+\right)  $. The heat kernel, as a
distribution, is given by
\[
Z_{\mathbb{R}}(x,t)=\mathcal{F}_{\xi\rightarrow x}^{-1}\left(  \exp
-t|\xi|_{\mathbb{R}}^{\gamma}\right)  ,\text{ for }x\in\mathbb{R}\text{,
}t>0.
\]
In order to have a probabilistic meaning, this kernel must be a probability
measure, then by a Bochner theorem, cf. \cite[Theorem 3.12]{Berg-Gunnar},
$\exp-t|\xi|_{\mathbb{R}}^{\gamma}$ is a positive definite function, and by a
theorem due to Schoenberg cf. \cite[Theorem 7.8]{Berg-Gunnar}, $|\xi
|_{\mathbb{R}}^{\gamma}$ is a negative definite function. Now, if
$\psi:\mathbb{R}\rightarrow\mathbb{C}$ is a negative definite function, then
$\left\vert \psi\left(  x\right)  \right\vert \leq C|x|_{\mathbb{R}}^{2}$ for
$|x|_{\mathbb{R}}\geq1$, cf. \cite[Corollary 7.16]{Berg-Gunnar}. This implies
that $0\leq\gamma\leq2$.

The family of `$p$-adic Laplacians' is very large, see e.g. \cite[Chapter
12]{KKZuniga}, \cite[Chapter 4]{Koch}, \cite[Chapter 3, Section XVI]{V-V-Z}
\cite[Chapter 2]{Zuniga-LNM-2016} and the references therein.

\section{The Model}

We fix $f$, $g$ $:\mathbb{R}^{2}\rightarrow\mathbb{R}$ two $\mathbb{R}%
$-analytic functions, and fix $d$, $\gamma$, $\alpha>0$. In this article we
consider the following non-Archimedean Turing system:%
\begin{equation}
\left\{
\begin{array}
[c]{l}%
u(\cdot,t),v(\cdot,t)\in L_{\mathbb{R}}^{2}\left(  B_{M}\right)
\cap\mathcal{L}\left(  B_{M}\right)  \text{, for }t\geq0\text{;}\\
\\
u(x,0),v(x,0)\in L_{\mathbb{R}}^{2}\left(  B_{M}\right)  \cap\mathcal{L}%
\left(  B_{M}\right)  \text{, }u(x,0),v(x,0)\geq0\text{;}\\
\\
\frac{\partial u}{\partial t}(x,t)=\gamma f\left(  u,v\right)  -\boldsymbol{D}%
_{x}^{\alpha}u(x,t)\text{;}\\
\\
\frac{\partial v}{\partial t}(x,t)=\gamma g\left(  u,v\right)
-d\boldsymbol{D}_{x}^{\alpha}v(x,t)\text{, }x\in B_{M}\text{, }t\geq0\text{.}%
\end{array}
\right.  \label{Turing_system}%
\end{equation}
Since we want to study self-organization patterns we need a `zero flux
boundary condition', which is
\begin{equation}
u(x,t),v(x,t)\equiv0\text{ for }x\in\mathbb{Q}_{p}\smallsetminus B_{M}\text{,
for any }t\geq0. \label{EQ_4}%
\end{equation}
Condition (\ref{EQ_4}) seems very strong in comparison with the classical one,
but there are at least to reasons for this choice. First, since $B_{M}$ is
open and closed, then its boundary is the empty set. Second, the operator
$\boldsymbol{D}_{x}^{\alpha}$ is non-local, as a consequence of this a
condition like $\boldsymbol{D}_{x}^{\alpha}u(x,t)\equiv0$ for any $x\in B_{M}%
$, for any $t\geq0$, it is not sufficient to stop the diffusion outside of
ball $B_{M}$, see (\ref{Eq_D_alpha}).

\section{Turing Instability Criteria}

We now consider a homogeneous steady state $\left(  u_{0},v_{0}\right)  $ of
(\ref{Turing_system}) which is a positive solution of
\begin{equation}
f(u,v)=g(u,v)=0. \label{EQ_6}%
\end{equation}
Since $u$, $v$ are real-valued functions, to study the linear stability of
$\left(  u_{0},v_{0}\right)  $ we can use the classical results, see e.g.
\cite[Chapter 2]{Murray-II}. Following Turing, in the absence of any spatial
variation, the homogeneous state must be linearly stable. With no spatial
variation $u$, $v$ satisfy%
\begin{equation}
\left\{
\begin{array}
[c]{c}%
\frac{\partial u}{\partial t}(x,t)=\gamma f\left(  u,v\right) \\
\\
\frac{\partial v}{\partial t}(x,t)=\gamma g\left(  u,v\right)  .
\end{array}
\right.  \label{EQ_7}%
\end{equation}
Notice that (\ref{EQ_7}) \ is an ordinary system of differential equations in
$\mathbb{R}^{2}$. In order to linearize about the steady state $\left(
u_{0},v_{0}\right)  $, we set%
\begin{equation}
\boldsymbol{w}=\left[
\begin{array}
[c]{c}%
w_{1}\\
\\
w_{2}%
\end{array}
\right]  =\left[
\begin{array}
[c]{c}%
u-u_{0}\\
\\
v-v_{0}.
\end{array}
\right]  \label{EQ_8}%
\end{equation}
By using the fact that $f$ and $g$ are $\mathbb{R}$-analytic, and assuming
that $\left\Vert \boldsymbol{w}\right\Vert _{L^{\infty}}:=\max\left\{
\left\Vert w_{1}\right\Vert _{L^{\infty}},\left\Vert w_{2}\right\Vert
_{L^{\infty}}\right\}  $ is small, then (\ref{EQ_7}) can be approximated as%
\begin{equation}
\frac{\partial\boldsymbol{w}}{\partial t}=\gamma\mathbb{J}\boldsymbol{w}%
\text{, } \label{EQ_9}%
\end{equation}
where \
\[
\mathbb{J}=\left[
\begin{array}
[c]{ccc}%
\frac{\partial f}{\partial u} &  & \frac{\partial f}{\partial v}\\
&  & \\
\frac{\partial g}{\partial u} &  & \frac{\partial g}{\partial v}%
\end{array}
\right]  \left(  u_{0},v_{0}\right)  =:\left[
\begin{array}
[c]{ccc}%
f_{u_{0}} &  & f_{v_{0}}\\
&  & \\
g_{u_{0}} &  & g_{v_{0}}%
\end{array}
\right]  .
\]
We now look for solutions of (\ref{EQ_9}) of the form%
\begin{equation}
\boldsymbol{w}\left(  t;\lambda\right)  =e^{\lambda t}\boldsymbol{w}_{0}.
\label{EQ_10}%
\end{equation}
By substituting (\ref{EQ_10}) in (\ref{EQ_9}), the eigenvalues $\lambda$ are
the solutions of
\[
\det\left(  \gamma\mathbb{J}-\lambda I\right)  =0,
\]
i.e.
\begin{equation}
\lambda^{2}-\gamma\left(  Tr\mathbb{J}\right)  \lambda+\gamma^{2}%
\det\mathbb{J}=0. \label{EQ_10AB}%
\end{equation}
Consequently%
\[
\lambda_{1,2}=\frac{\gamma}{2}\left\{  Tr\mathbb{J}\text{ }\pm\sqrt{\left(
Tr\mathbb{J}\right)  ^{2}-4\det\mathbb{J}}\right\}  .
\]
The steady state $\boldsymbol{w}=\boldsymbol{0}$ is linearly stable if
$\operatorname{Re}\lambda_{1,2}<0$, this last condition is guaranteed if
\begin{equation}
Tr\mathbb{J}<0\text{ \ and \ }\det\mathbb{J}>0. \label{EQ_10A}%
\end{equation}
We linearize the full reaction-ultradiffusion system about the steady state,
which is $\boldsymbol{w}=\boldsymbol{0}:\boldsymbol{=}\left[
\begin{array}
[c]{c}%
0\\
0
\end{array}
\right]  $, see (\ref{EQ_8}), to get%
\begin{equation}
\frac{\partial\boldsymbol{w}}{\partial t}(x,t)=\gamma\mathbb{J}\boldsymbol{w}%
(x,t)\boldsymbol{-}{\scriptsize D}\boldsymbol{D}_{x}^{\alpha}\boldsymbol{w}%
(x,t), \label{EQ_11}%
\end{equation}
where
\begin{equation}
{\scriptsize D=}\left[
\begin{array}
[c]{ccc}%
1 &  & 0\\
&  & \\
0 &  & d
\end{array}
\right]  \text{, \ }\boldsymbol{D}_{x}^{\alpha}\boldsymbol{w}:\boldsymbol{=}%
\left[
\begin{array}
[c]{c}%
\boldsymbol{D}_{x}^{\alpha}\boldsymbol{w}_{1}\\
\\
\boldsymbol{D}_{x}^{\alpha}\boldsymbol{w}_{2}%
\end{array}
\right]  \text{.} \label{EQ_13}%
\end{equation}
To solve the system (\ref{EQ_11}) subject to the boundary conditions
(\ref{EQ_4}), we first determine a solution $\boldsymbol{w}_{\kappa}$ of the
following eigenvalue problem:%
\begin{equation}
\left\{
\begin{array}
[c]{l}%
\boldsymbol{D}_{x}^{\alpha}\boldsymbol{w}_{\kappa}(x)=\kappa\boldsymbol{w}%
_{\kappa}(x)\\
\\
\boldsymbol{w}_{\kappa}\in L_{\mathbb{R}}^{2}\left(  B_{M}\right)
\cap\mathcal{L}\left(  B_{M}\right)  .
\end{array}
\right.  \label{EQ_14}%
\end{equation}
In Section \ref{Section_Eigenvalue_problems}, the existence of solution for
the eigenvalue problem (\ref{EQ_14}) was established. Indeed, $\ $if
$\boldsymbol{w}_{\kappa}=\left[
\begin{array}
[c]{l}%
w_{1,\kappa}\\
w_{2,\kappa}%
\end{array}
\right]  $, then
\begin{gather*}
w_{1,\kappa},w_{2,\kappa}\in%
{\textstyle\bigsqcup_{rnj}}
\left\{  p^{\frac{-r}{2}}\cos\left(  \left\{  p^{-1}j\left(  p^{r}x-n\right)
\right\}  _{p}\right)  \Omega\left(  \left\vert p^{r}x-n\right\vert
_{p}\right)  \right\}
{\textstyle\bigsqcup}
\\%
{\textstyle\bigsqcup_{rnj}}
\left\{  p^{\frac{-r}{2}}\sin\left(  \left\{  p^{-1}j\left(  p^{r}x-n\right)
\right\}  _{p}\right)  \Omega\left(  \left\vert p^{r}x-n\right\vert
_{p}\right)  \right\}
{\textstyle\bigsqcup}
\left\{  \Omega\left(  p^{-M}\left\vert x\right\vert _{p}\right)  \right\}  ,
\end{gather*}
with $r=r(\kappa)$, $j\in\left\{  1,\cdots,p-1\right\}  $, and $n$ $\in$
$\mathbb{Q}_{p}/\mathbb{Z}_{p}$ as in (\ref{conditions}).The value
$\frac{p^{-\alpha M}}{1-p^{-\alpha-1}}$ is the eigenvalue corresponding to
$\Omega\left(  p^{-M}\left\vert x\right\vert _{p}\right)  $. We now look for a
solution $\boldsymbol{w}(x,t)$ of (\ref{EQ_11}) of the form $\boldsymbol{w}%
(x,t)=\sum_{\kappa,\lambda}C_{\kappa,\lambda}e^{\lambda t}\boldsymbol{w}%
_{\kappa}\left(  x\right)  $. The function $e^{\lambda t}\boldsymbol{w}%
_{\kappa}\left(  x\right)  $ is a non-trivial solution of (\ref{EQ_11}) if
$\lambda$ satisfies%
\begin{equation}
\det\left(  \lambda I-\gamma\mathbb{J}+\kappa{\scriptsize D}\right)  =0,
\label{EQ_16A}%
\end{equation}
i.e.
\begin{equation}
\lambda^{2}+\left\{  \kappa\left(  1+d\right)  -\gamma Tr\mathbb{J}\right\}
\lambda+h\left(  \kappa\right)  =0, \label{EQ_16}%
\end{equation}
where%
\begin{equation}
h\left(  \kappa\right)  :=d\kappa^{2}-\gamma\kappa\left(  df_{u_{0}}+g_{v_{0}%
}\right)  +\gamma^{2}\det\mathbb{J}. \label{EQ_17}%
\end{equation}
In the Archimedean case, see e.g. \cite[Section 2.3]{Murray-II}, condition
(\ref{EQ_16A}) becomes condition (\ref{EQ_10AB}) when $\kappa=0$. Since
$\kappa=0$ is not an eigenvalue of operator $\boldsymbol{D}_{x}^{\alpha}$,
conditions (\ref{EQ_16A}) and (\ref{EQ_10AB}) are independent.

The steady state $\left(  u_{0},v_{0}\right)  $ is linearly stable if both
solutions of (\ref{EQ_16}) have $\operatorname{Re}\left(  \lambda\right)  <0$.
Conditions (\ref{EQ_10A}) guarantee that the steady state is stable in absence
of spatial effects, i.e. $\operatorname{Re}\left(  \lambda\mid_{\kappa
=0}\right)  <0$. For the steady state to be unstable to spatial disturbances
we require $\operatorname{Re}\left(  \lambda\left(  \kappa\right)  \right)
>0$ for some $\kappa\neq0$. This can happen if either the coefficient of
$\lambda$ in (\ref{EQ_16}) is negative, or if $h\left(  \kappa\right)  <0$ for
some $\kappa\neq0$ in (\ref{EQ_17}). Since $Tr\mathbb{J}<0$ from conditions
(\ref{EQ_10A}) and the coefficient of $\lambda$ in (\ref{EQ_17}) is $\left(
1+d\right)  -\gamma Tr\mathbb{J}$, which is positive, therefore, the only way
$\operatorname{Re}\left(  \lambda\left(  \kappa\right)  \right)  $ can be
positive is if $h\left(  \kappa\right)  <0$ for some $\kappa\neq0$. Since
$\det\mathbb{J}>0$ from (\ref{EQ_10A}), in order $h\left(  \kappa\right)  $ to
be \ negative, it is necessary that $\left(  df_{u_{0}}+g_{v_{0}}\right)  >0$.
Now, since $f_{u_{0}}+g_{v_{0}}=Tr\mathbb{J}<0$, necessarily $d\neq1$ and
$f_{u_{0}}$ and $g_{v_{0}}$ must have opposite signs. So an additional
requirement to those of (\ref{EQ_10A}) is
\begin{equation}
d\neq1. \label{EQ_18}%
\end{equation}
This is a necessary condition, but not sufficient for $\operatorname{Re}%
\left(  \lambda\left(  \kappa\right)  \right)  >0$. For $h\left(
\kappa\right)  $ to be negative for some nonzero $\kappa$, the minimum
$h_{\text{min}}$ of $h\left(  \kappa\right)  $ must be negative. An elementary
calculation shows that
\begin{equation}
h_{\text{min}}=\gamma^{2}\left\{  \det\mathbb{J}-\frac{\left(  df_{u_{0}%
}+g_{v_{0}}\right)  ^{2}}{4d}\right\}  , \label{EQ_19}%
\end{equation}
and the minimum is achieved at
\begin{equation}
\kappa_{\text{min}}=\gamma\frac{df_{u_{0}}+g_{v_{0}}}{2d} \label{EQ_20}%
\end{equation}
Thus the condition $h\left(  \kappa\right)  <0$ for some $\kappa\neq0$ is
\begin{equation}
\frac{\left(  df_{u_{0}}+g_{v_{0}}\right)  ^{2}}{4d}>\det\mathbb{J}.
\label{EQ_21}%
\end{equation}

A bifurcation occurs when $h_{\text{min}}=0$, see (\ref{EQ_19}), for fixed
kinetics parameters, this condition,%
\begin{equation}
\det\mathbb{J}=\frac{\left(  df_{u_{0}}+g_{v_{0}}\right)  ^{2}}{4d},
\label{EQ_26A}%
\end{equation}
defines a critical diffusion $d_{c}$, which is given as an appropriate root
of
\begin{equation}
f_{u_{0}}^{2}d_{c}^{2}+2\left(  2f_{v_{0}}g_{u_{0}}-f_{u_{0}}g_{v_{0}}\right)
d_{c}+g_{v_{0}}^{2}=0. \label{EQ_26B}%
\end{equation}

The model for $d>d_{c}$ exhibits Turing instability, while for $d<d_{c}$ no.
Notice that $d_{c}>1$. A critical `wavenumber' $\kappa_{c}$ is obtained by
using (\ref{EQ_20}):%
\begin{equation}
\kappa_{c}=\gamma\frac{d_{c}f_{u_{0}}+g_{v_{0}}}{2d_{c}}=\gamma\sqrt
{\frac{\det\mathbb{J}}{d_{c}}}. \label{EQ_26C}%
\end{equation}
When $d>d_{c}$, there exists a range of unstable of positive wavenumbers
\ $\kappa_{1}<\kappa<\kappa_{2}$, where $\kappa_{1}$, $\kappa_{2}$ are the
zeros of $h\left(  \kappa\right)  =0$, see (\ref{EQ_17}):
\begin{equation}
\kappa_{1}=\frac{\gamma}{2d}\left\{  \left(  df_{u_{0}}+g_{v_{0}}\right)
-\sqrt{\left(  df_{u_{0}}+g_{v_{0}}\right)  ^{2}-4d\det\mathbb{J}}\right\}  ,
\label{Kappa_1}%
\end{equation}%
\begin{equation}
\kappa_{2}=\frac{\gamma}{2d}\left\{  \left(  df_{u_{0}}+g_{v_{0}}\right)
+\sqrt{\left(  df_{u_{0}}+g_{v_{0}}\right)  ^{2}-4d\det\mathbb{J}}\right\}  .
\label{Kappa_2}%
\end{equation}
We call the function $\lambda\left(  \kappa\right)  $ the \textit{dispersion
relation}. Notice that, within the unstable range, $\operatorname{Re}%
\lambda\left(  \kappa\right)  >0$ has a maximum for the wavenumber
$\kappa_{\text{min}}^{\left(  0\right)  }$ obtained from (\ref{EQ_20}) with
$d>d_{c}$. Some typical plots for $\lambda\left(  \kappa\right)  $ and
$\operatorname{Re}\lambda\left(  \kappa\right)  $ are showed in \cite[Section
2.3]{Murray-II}, see also the figure 2.5 in \cite[Section 2.3]{Murray-II}.
Then as $t$ increases the behavior of $\boldsymbol{w}\left(  x,t\right)  $ is
controlled by the dominant modes, i.e. those $e^{\lambda\left(  \kappa\right)
t}\boldsymbol{w}_{\kappa}(x)$ with $\operatorname{Re}\lambda\left(
\kappa\right)  >0$, since the other modes tend to zero exponentially. Then
\begin{gather}
\boldsymbol{w}\left(  x,t\right)  \sim\text{ }\sum\limits_{\kappa_{1}%
<\kappa<\kappa_{2}}\sum\limits_{nj}\boldsymbol{A}_{rnj}e^{\lambda t}%
p^{\frac{-r}{2}}\cos\left(  \left\{  p^{-1}j\left(  p^{r}x-n\right)  \right\}
_{p}\right)  \Omega\left(  \left\vert p^{r}x-n\right\vert _{p}\right)
+\nonumber\\
\sum\limits_{\kappa_{1}<\kappa<\kappa_{2}}\sum\limits_{nj}\boldsymbol{B}%
_{rnj}e^{\lambda t}p^{\frac{-r}{2}}\sin\left(  \left\{  p^{-1}j\left(
p^{r}x-n\right)  \right\}  _{p}\right)  \Omega\left(  \left\vert
p^{r}x-n\right\vert _{p}\right)  \text{ } \label{EQ_Expansion}%
\end{gather}
for $t\rightarrow+\infty$. In the above expansion in all the series the $r$s
and $j$s take only a finite numbers of values. Indeed, all the $\kappa$s but
one has the form $p^{\left(  1-r\right)  \alpha}$, the condition $\kappa
_{1}<\kappa<\kappa_{2}$ implies that there is only a finite number of $r$s.
But the $n$s run through an infinite set, see (\ref{conditions}). We now fix
$r$, then, for a given $x\in B_{M}$ there exists only a finite numbers of
balls of type $B_{r}\left(  p^{-r}n\right)  $ containing $x$, see
(\ref{conditions}). This fact implies that the value at $\left(  x,t\right)  $
of $\boldsymbol{w}\left(  x,t\right)  $ in expansion (\ref{EQ_Expansion}) is
determined only by a finite number of $n$s, and consequently the series in
expansion (\ref{EQ_Expansion}) is convergent. We now formulate the Turing
instability criteria for our model.

\begin{notation}
We denote by $\sigma\left(  \boldsymbol{D}^{\alpha}\right)  $\ the spectra of
$\boldsymbol{D}^{\alpha}$.
\end{notation}

\begin{theorem}
\label{Theorem1}Consider the reaction-diffusion system (\ref{Turing_system}).
The steady state $(u_{0},v_{0})$ is linearly unstable (Turing unstable) if the
following conditions hold:

\noindent(T1) $\mathit{Tr}A\mathbb{=}f_{u_{0}}+g_{v_{0}}<0$;

\noindent(T2) $\det\mathbb{A}=f_{u_{0}}g_{v_{0}}-f_{v_{0}}g_{u_{0}}>0$;

\noindent(T3) $df_{u_{0}}+g_{v_{0}}>0$;

\noindent(T4) the derivatives $f_{u_{0}}$ and $g_{v_{0}}$ must have opposite signs;

\noindent(T5) $\left(  df_{u_{0}}+g_{v_{0}}\right)  ^{2}-4d\left(  f_{u_{0}%
}g_{v_{0}}-f_{v_{0}}g_{u_{0}}\right)  >0$;

\noindent(T6) $\left\{  \kappa\in\sigma\left(  \boldsymbol{D}^{\alpha}\right)
;\kappa_{1}<\kappa<\kappa_{2}\right\}  \neq\emptyset$.

Moreover, there are infinitely many unstable eigenmodes, and the Turing
pattern has the form (\ref{EQ_Expansion}).
\end{theorem}

\section{\label{Section_Discrete_models}Discrete Models of Morphogenesis}

\subsection{The spaces $\mathcal{D}_{M}^{-L}$}

We fix $M\in\mathbb{Z}$ and $L\in\mathbb{N}$, with $L\geq-M$, and define%
\[
G_{L,M}=p^{-M}\mathbb{Z}_{p}/p^{L}\mathbb{Z}_{p}.
\]
Then, $G_{L,M}$ is a finite ring, with $\#G_{L,M}=p^{L+M}$ elements. We fix
the following a set of representatives for $G_{L,M}$:%
\[
I=I_{-M}p^{-M}+I_{-M+1}p^{-M+1}+\ldots+I_{L-1}p^{L-1},
\]
where the $I_{j}$s are $p$-adic digits, i.e. elements from $\left\{
0,1,\ldots,p-1\right\}  $.

On the other hand, since $I=p^{-M}\widetilde{I}$, with
\[
\widetilde{I}=I_{-M}+I_{-M+1}p+\ldots+I_{L-1}p^{L+M-1}\in\mathbb{Z}%
_{p}/p^{L+M}\mathbb{Z}_{p},
\]
and elements of $\mathbb{Z}_{p}/p^{L+M}\mathbb{Z}_{p}$ form in a natural way a
rooted tree with $L+M$-levels, where the level $1\leq j\leq L+M$ is formed by
the elements of the form $a_{0}+\ldots+a_{j-1}p^{j-1}$. \ In this way we can
identify $G_{L,M}$ with a rooted tree. This a identification is not unique,
see e.g. \cite{zuniga-Nonlieal}.

By considering $G_{L,M}$ as a subset of $\mathbb{Q}_{p}$, we can endow
$G_{L,M}$ with a norm denoted as $\left\vert \cdot\right\vert _{p}$, and thus
$\left(  G_{L,M},\left\vert \cdot\right\vert _{p}\right)  $ becomes a finite
ultrametric space.

We define $\mathcal{D}_{M}^{-L}$ to be the $\mathbb{R}$-vector space formed by
test functions $\varphi$ supported in the ball $p^{-M}\mathbb{Z}_{p}$ having
the form
\[
\varphi\left(  x\right)  =p^{\frac{L}{2}}\sum\limits_{I\in G_{L,M}}%
\varphi\left(  I\right)  \Omega\left(  p^{L}\left\vert x-I\right\vert
_{p}\right)  \text{, with }\varphi\left(  I\right)  \in\mathbb{R}.
\]
Since $\Omega\left(  p^{L}\left\vert x-I\right\vert _{p}\right)  \Omega\left(
p^{L}\left\vert x-J\right\vert _{p}\right)  =0$ if $I\neq J$, the set
\[
\left\{  p^{\frac{L}{2}}\Omega\left(  p^{L}\left\vert x-I\right\vert
_{p}\right)  ;I\in G_{L,M}\right\}
\]
is an orthonormal basis for $\mathcal{D}_{M}^{-L}$. Now, by using that%
\begin{align*}
\left\Vert \varphi\right\Vert _{L^{2}}  &  =\sqrt{p^{L}\sum\limits_{I\in
G_{L,M}}\left\vert \varphi\left(  I\right)  \right\vert ^{2}\int
\limits_{p^{-M}\mathbb{Z}_{p}}\Omega\left(  p^{L}\left\vert x-I\right\vert
_{p}\right)  dx}\\
&  =\sqrt{\sum\limits_{I\in G_{L,M}}\left\vert \varphi\left(  I\right)
\right\vert ^{2},}%
\end{align*}
we have
\[
\left(  \mathcal{D}_{M}^{-L},\left\Vert \cdot\right\Vert _{L^{2}}\right)
\simeq\left(  \mathbb{R}^{\#G_{L,M}},\left\Vert \cdot\right\Vert _{\mathbb{R}%
}\right)  \text{ as Banach spaces,}%
\]
where $\left\Vert \cdot\right\Vert _{\mathbb{R}}$ denotes the usual norm of
$\mathbb{R}^{\#G_{L,M}}$.

\subsection{Discretization of the operator $\boldsymbol{D}^{\alpha}$}

We set $\lambda_{M}:=\frac{\left(  1-p^{-1}\right)  p^{-\alpha M}%
}{1-p^{-\alpha-1}}>0$. Notice that
\[
\boldsymbol{D}_{x}^{\alpha}\varphi(x)=\left(  \boldsymbol{D}_{M}^{\alpha
}-\lambda_{M}\right)  \varphi(x)\text{ \ for }\varphi\in\mathcal{D}%
_{\mathbb{R}}(p^{-M}\mathbb{Z}_{p})\text{,}%
\]
where%
\[
\boldsymbol{D}_{M}^{\alpha}\varphi(x):=\frac{1-p^{\alpha}}{1-p^{-\alpha-1}%
}\int\limits_{p^{-M}\mathbb{Z}_{p}}\frac{\varphi\left(  x-y\right)
-\varphi\left(  x\right)  }{\left\vert y\right\vert _{p}^{\alpha+1}}dy\text{.}%
\]
The operator $\boldsymbol{D}_{M}^{\alpha}$ $-\lambda_{M}$ is a non-negative,
symmetric operator on $L_{\mathbb{R}}^{2}(p^{-M}\mathbb{Z}_{p}) $.
Furthermore, its closure, also denoted by $\boldsymbol{D}_{M}^{\alpha}%
$,$-\lambda_{M}$ is a self-adjoint operator, see \cite[Section 3.3.2]{Koch}.

Every wavelet $\Psi_{rnj}\left(  x\right)  $, with support in $p^{-M}%
\mathbb{Z}_{p}$, with $r$, $n$, $j$ satisfying (\ref{conditions}), is an
eigenfunction of $\boldsymbol{D}_{M}^{\alpha}$ with eigenvalue $p^{\left(
1-r\right)  \alpha}$. In addition $p^{\frac{M}{2}}\Omega\left(  p^{-M}%
\left\vert x\right\vert _{p}\right)  $ is also a eigenfunction of
$\boldsymbol{D}_{M}^{\alpha}$ with eigenvalue $\lambda_{M}p^{\frac{M}{2}}$.

As a discretization of $\boldsymbol{D}_{M}^{\alpha}$ $-\lambda_{M}$, we pick
its restriction to $\mathcal{D}_{M}^{-L}$, which \ is denoted as
$\boldsymbol{D}_{L,M}^{\alpha}-\lambda_{M}$. \ Since $\mathcal{D}_{M}^{-L}$ is
a finite vector space $\boldsymbol{D}_{L,M}^{\alpha}-\lambda_{M}$ is
represented by a matrix $A_{L,M}^{\alpha}$.

\subsection{Computation of the matrix $A_{L,M}^{\alpha}$}

In order to compute the matrix $A_{L,M}^{\alpha}$, we first compute%
\begin{align*}
&  \boldsymbol{D}_{M}^{\alpha}\left(  p^{\frac{L}{2}}\Omega\left(
p^{L}\left\vert x-I\right\vert _{p}\right)  \right) \\
&  =p^{\frac{L}{2}}\frac{1-p^{\alpha}}{1-p^{-\alpha-1}}\int\limits_{p^{-M}%
\mathbb{Z}_{p}}\frac{\Omega\left(  p^{L}\left\vert x-y-I\right\vert
_{p}\right)  -\Omega\left(  p^{L}\left\vert x-I\right\vert _{p}\right)
}{\left\vert y\right\vert _{p}^{\alpha+1}}dy\\
&  =p^{\frac{L}{2}}\frac{1-p^{\alpha}}{1-p^{-\alpha-1}}\sum\limits_{J\in
G_{L,M}}\text{ }\int\limits_{\left\vert y-J\right\vert _{p}\leq p^{-L}}%
\frac{\Omega\left(  p^{L}\left\vert x-y-I\right\vert _{p}\right)
-\Omega\left(  p^{L}\left\vert x-I\right\vert _{p}\right)  }{\left\vert
y\right\vert _{p}^{\alpha+1}}dy\\
&  =:p^{\frac{L}{2}}\frac{1-p^{\alpha}}{1-p^{-\alpha-1}}\sum\limits_{J\in
G_{L,M}}\mathcal{I}_{J}\left(  x,I,L\right)  \text{.}%
\end{align*}
We now compute the integrals $\mathcal{I}_{J}\left(  x,I,L\right)  $. We
consider first the case $J\neq0$. By using that
\[
\Omega\left(  p^{L}\left\vert x-I\right\vert _{p}\right)  \ast\Omega\left(
p^{L}\left\vert x-J\right\vert _{p}\right)  =p^{-L}\Omega\left(
p^{L}\left\vert x-\left(  I+J\right)  \right\vert _{p}\right)  ,
\]
we have%
\begin{gather*}
\mathcal{I}_{J}\left(  x,I,L\right)  =\frac{1}{\left\vert J\right\vert
_{p}^{\alpha+1}}\left\{  \Omega\left(  p^{L}\left\vert x-I\right\vert
_{p}\right)  \ast\Omega\left(  p^{L}\left\vert x-J\right\vert _{p}\right)
-p^{-L}\Omega\left(  p^{L}\left\vert x-I\right\vert _{p}\right)  \right\} \\
=\frac{p^{-L}}{\left\vert J\right\vert _{p}^{\alpha+1}}\left\{  \Omega\left(
p^{L}\left\vert x-\left(  I+J\right)  \right\vert _{p}\right)  -\Omega\left(
p^{L}\left\vert x-I\right\vert _{p}\right)  \right\}  .
\end{gather*}

Now in the case $J=0$, we have $\mathcal{I}_{J}\left(  x,I,L\right)  =0$.
Indeed,%
\[
\mathcal{I}_{0}\left(  x,I,L\right)  =\int\limits_{\left\vert y\right\vert
_{p}\leq p^{-L}}\frac{\Omega\left(  p^{L}\left\vert x-y-I\right\vert
_{p}\right)  -\Omega\left(  p^{L}\left\vert x-I\right\vert _{p}\right)
}{\left\vert y\right\vert _{p}^{\alpha+1}}dy.
\]
If $\Omega\left(  p^{L}\left\vert x-I\right\vert _{p}\right)  =1$, then $x\in
I+p^{L}\mathbb{Z}_{p}$ and since $y\in p^{L}\mathbb{Z}_{p}$ we have $x-y-I\in
p^{L}\mathbb{Z}_{p}$, which implies that $\Omega\left(  p^{L}\left\vert
x-y-I\right\vert _{p}\right)  =1$. Now, if $\Omega\left(  p^{L}\left\vert
x-I\right\vert _{p}\right)  =0$, i.e. if $x\in K+p^{L}\mathbb{Z}_{p}$ for some
$K\neq I$, then by using that $p^{L}\mathbb{Z}_{p}\cap K-I+p^{L}\mathbb{Z}%
_{p}=\varnothing$, we have
\[
\mathcal{I}_{0}\left(  x,I,L\right)  =\int\limits_{p^{L}\mathbb{Z}_{p}\cap
K-I+p^{L}\mathbb{Z}_{p}}\frac{1}{\left\vert y\right\vert _{p}^{\alpha+1}%
}dy=0.
\]
Now, we use the fact that $G_{L,M}$ is a additive group to conclude that
\begin{gather*}
\boldsymbol{D}_{M}^{\alpha}\left(  p^{\frac{L}{2}}\Omega\left(  p^{L}%
\left\vert x-I\right\vert _{p}\right)  \right) \\
=p^{-\frac{L}{2}}\frac{1-p^{\alpha}}{1-p^{-\alpha-1}}\sum
\limits_{\substack{K\in G_{L,M} \\K\neq I}}\text{ }\frac{1}{\left\vert
K-I\right\vert _{p}^{\alpha+1}}\Omega\left(  p^{L}\left\vert x-K\right\vert
_{p}\right) \\
-p^{-\frac{L}{2}}\frac{1-p^{\alpha}}{1-p^{-\alpha-1}}\left(  \sum\limits_{
_{\substack{K\in G_{L,M} \\K\neq I}}}\text{ }\frac{1}{\left\vert
K-I\right\vert _{p}^{\alpha+1}}\right)  \Omega\left(  p^{L}\left\vert
x-I\right\vert _{p}\right)  ,
\end{gather*}
and the entries of the matrix $A_{L,M}^{\alpha}=\left[  A_{K,I}^{\alpha
}\right]  _{_{K,I\in G_{L,M}}}$ are given as%
\begin{equation}
A_{K,I}^{\alpha}=\left\{
\begin{array}
[c]{ccc}%
p^{-\frac{L}{2}}\frac{1-p^{\alpha}}{1-p^{-\alpha-1}}\text{ }\frac
{1}{\left\vert K-I\right\vert _{p}^{\alpha+1}} & \text{if} & K\neq I\\
&  & \\
-p^{-\frac{L}{2}}\frac{1-p^{\alpha}}{1-p^{-\alpha-1}}\sum\limits_{K\neq
I}\frac{1}{\left\vert K-I\right\vert _{p}^{\alpha+1}}-\lambda_{M} & \text{if}
& K=I.
\end{array}
\right.  \label{EQ_A_matrix}%
\end{equation}

\subsection{Discretization of the $p$-adic Turing System}

In the discretization of the Turing system (\ref{Turing_system}), we use the
following approximation for functions $u(x,t)$, $v(x,t)$:%
\begin{equation}
u^{\left(  L\right)  }(x,t)=\sum\limits_{I\in G_{L,M}}u^{\left(  L\right)
}\left(  I,t\right)  \Omega\left(  p^{L}\left\vert x-I\right\vert _{p}\right)
\label{Eq_dis_1}%
\end{equation}
and
\begin{equation}
v^{\left(  L\right)  }(x,t)=\sum\limits_{I\in G_{L,M}}v^{\left(  L\right)
}\left(  I,t\right)  \Omega\left(  p^{L}\left\vert x-I\right\vert _{p}\right)
, \label{Eq_dis_2}%
\end{equation}
where $u^{\left(  L\right)  }\left(  I,\cdot\right)  $, $v^{\left(  L\right)
}\left(  I,\cdot\right)  \in C^{1}\left(  \left[  0,T\right]  \right)  $ for
some fixed positive $T$. Furthermore, we set
\[
u^{\left(  L\right)  }(x,t)=\left[  u^{\left(  L\right)  }\left(  I,t\right)
\right]  _{I\in G_{L,M}}\text{, \ }v^{\left(  L\right)  }(x,t)=\left[
v^{\left(  L\right)  }\left(  I,t\right)  \right]  _{I\in G_{L,M}}.
\]
We assume that range$\left(  u^{\left(  L\right)  }(x,t)\right)  \times$range
$\left(  v^{\left(  L\right)  }(x,t)\right)  $ is contained in the domains of
convergence of $f$, $g$. Then%
\begin{gather*}
f\left(  \sum\limits_{I\in G_{L,M}}u^{\left(  L\right)  }\left(  I,t\right)
\Omega\left(  p^{L}\left\vert x-I\right\vert _{p}\right)  ,\sum\limits_{J\in
G_{L,M}}v^{\left(  L\right)  }\left(  J,t\right)  \Omega\left(  p^{L}%
\left\vert x-J\right\vert _{p}\right)  \right)  =\\
\sum\limits_{I\in G_{L,M}}f\left(  u^{\left(  L\right)  }\left(  I,t\right)
,v^{\left(  L\right)  }\left(  I,t\right)  \right)  \Omega\left(
p^{L}\left\vert x-I\right\vert _{p}\right)  .
\end{gather*}
A similar formula holds for function $g$. Then the discretization of the
$p$-adic Turing system has the form:%
\begin{multline*}
\frac{\partial}{\partial t}\left[  u^{\left(  L\right)  }\left(  I,t\right)
\right]  _{I\in G_{L,M}}=\gamma\left[  f\left(  u^{\left(  L\right)  }\left(
I,t\right)  ,v^{\left(  L\right)  }\left(  I,t\right)  \right)  \right]
_{_{I\in G_{L,M}}}\text{ }\\
-A_{L,M}^{\alpha}\left[  u^{\left(  L\right)  }\left(  I,t\right)  \right]
_{I\in G_{L,M}}%
\end{multline*}%
\begin{multline*}
\frac{\partial}{\partial t}\left[  v^{\left(  L\right)  }\left(  I,t\right)
\right]  _{I\in G_{L,M}}=\gamma\left[  g\left(  u^{\left(  L\right)  }\left(
I,t\right)  ,v^{\left(  L\right)  }\left(  I,t\right)  \right)  \right]
_{_{I}\in G_{L,M}}\\
-dA_{L,M}^{\alpha}\left[  v^{\left(  L\right)  }\left(  I,t\right)  \right]
_{I\in G_{L,M}},
\end{multline*}
where $I\in G_{L,M}$. Equivalently,
\begin{gather*}
\frac{\partial}{\partial t}u^{\left(  L\right)  }\left(  I,t\right)  =\gamma
f\left(  u^{\left(  L\right)  }\left(  I,t\right)  ,v^{\left(  L\right)
}\left(  I,t\right)  \right)  -p^{-\frac{L}{2}}\frac{1-p^{\alpha}%
}{1-p^{-\alpha-1}}\sum\limits_{J\neq I}\text{ }\frac{u^{\left(  L\right)
}\left(  J,t\right)  }{\left\vert J-I\right\vert _{p}^{\alpha+1}}\\
+p^{-\frac{L}{2}}\frac{1-p^{\alpha}}{1-p^{-\alpha-1}}\left(  \sum
\limits_{_{J\neq I}}\text{ }\frac{1}{\left\vert J-I\right\vert _{p}^{\alpha
+1}}-\frac{\left(  1-p^{-1}\right)  p^{-\alpha M+\frac{L}{2}}}{1-p^{-\alpha}%
}\right)  u^{\left(  L\right)  }\left(  I,t\right)  ,
\end{gather*}%
\begin{gather*}
\frac{\partial}{\partial t}v^{\left(  L\right)  }\left(  I,t\right)  =\gamma
g\left(  u^{\left(  L\right)  }\left(  I,t\right)  ,v^{\left(  L\right)
}\left(  I,t\right)  \right)  -dp^{-\frac{L}{2}}\frac{1-p^{\alpha}%
}{1-p^{-\alpha-1}}\sum\limits_{J\neq I}\text{ }\frac{u^{\left(  L\right)
}\left(  J,t\right)  }{\left\vert J-I\right\vert _{p}^{\alpha+1}}\\
+dp^{-\frac{L}{2}}\frac{1-p^{\alpha}}{1-p^{-\alpha-1}}\left(  \sum
\limits_{_{J\neq I}}\text{ }\frac{1}{\left\vert J-I\right\vert _{p}^{\alpha
+1}}-\frac{\left(  1-p^{-1}\right)  p^{-\alpha M+\frac{L}{2}}}{1-p^{-\alpha}%
}\right)  u^{\left(  L\right)  }\left(  I,t\right)  ,
\end{gather*}
where $I\in G_{L,M}$.

\end{document}